\documentclass[10pt,a4paper]{siamltex}
\pdfoutput=1
\usepackage{amsmath,amssymb}
\usepackage
   {graphicx}
\title{On convex functions and the finite element method}
\author{N\'estor E. Aguilera\footnotemark[2] \and Pedro Morin\footnotemark[3]%
      }

\usepackage[%
   pdfpagemode=UseNone,
   bookmarks=true,
   pagebackref=true, 
   colorlinks,
   linkcolor=blue,
   anchorcolor=blue,
   citecolor=blue,
   filecolor=blue,
   pagecolor=blue,
   urlcolor=blue
   ]{hyperref}
\renewcommand*{\backref}[1]{}
\renewcommand*{\backrefalt}[4]{%
   \ifcase #1 %
      {} 
   \or
      (\textit{Page~#2}) 
   \else
      (\textit{Pages~#2}) 
   \fi}
\newcommand{\FE}{FE}

\newcommand{\strutmio}[2][0pt]{
   \rule[#1]{0pt}{#2}}
  \newtheorem{assumptions}[theorem]{Assumptions}
  \newtheorem{example}[theorem]{Example}
\newcommand{\ENDsymbol}{$\Diamond$}
\newcommand{\END}{~\hfill\ENDsymbol}

\newcommand{\RR}{\mathbb R}
\newcommand{\boundary}{\partial}
\newcommand{\clausure}[1]{\overline{#1}}

\newcommand{\Conv}{\mathcal{C}}

\newcommand{\Mesh}[1]{\mathcal{M}^{#1}}
\newcommand{\Triang}[1]{\mathcal{T}^h}
\newcommand{\IndV}[1]{{I^{#1}_{\text{trial}}}}
\newcommand{\IndW}[1]{{I^{#1}_{\text{test}}}}

\newcommand{\downto}{\ensuremath{\downarrow}}
\newcommand{\dd}{\mathrm{d}}
\DeclareMathOperator{\dist}{dist}
\newcommand{\interp}[1]{\mathcal{I}^{#1}} 
\newcommand{\abs}[1]{{\lvert#1\rvert}}
\newcommand{\norm}[1]{{\lVert#1\rVert}}

\newcommand{\degree}{^\circ}
\newcommand{\intern}[1]{\langle#1\rangle}
\newcounter{savectr}

\newcommand{\aut}[1]{\textsc{#1}}
\newcommand{\tita}[1]{#1}           
\newcommand{\titb}[1]{\emph{#1}}    
\newcommand{\titj}[1]{\emph{#1}}    
\usepackage{xspace}

\begin{document}
\maketitle

\renewcommand{\thefootnote}{\fnsymbol{footnote}}
\footnotetext[2]{Consejo Nacional de Investigaciones Cient\'{\i}ficas y T\'{e}cnicas and Universidad Nacional del Litoral, Argentina, \url{aguilera@santafe-conicet.gov.ar}. Partially supported by CONICET (Argentina) Grant PIP 5810.}
\footnotetext[3]{Consejo Nacional de Investigaciones Cient\'{\i}ficas y T\'{e}cnicas and Universidad Nacional del Litoral, Argentina, \url{pmorin@santafe-conicet.gov.ar}. Partially supported by CONICET (Argentina) through Grant PIP 5478, and Universidad Nacional del Litoral through Grant CAI+D 12/H421}

\begin{abstract}
Many problems of theoretical and practical interest involve finding a convex or concave function.
For instance, optimization problems such as finding the projection on the convex functions in $H^k(\Omega)$, or some problems in economics.

In the continuous setting and assuming smoothness, the convexity constraints may be given locally by asking the Hessian matrix to be positive semidefinite, but in making discrete approximations two difficulties arise: the continuous solutions may be not smooth, and an adequate discrete version of the Hessian must be given.

In this paper we propose a finite element description of the Hessian, and prove convergence under very general conditions, even when the continuous solution is not smooth, working on any dimension, and requiring a linear number of constraints in the number of nodes.

Using semidefinite programming codes, we show concrete examples of approximations to optimization problems.

\end{abstract}


\section{Introduction}
\label{sec:intro}

Convex and concave functions appear naturally in many disciplines of science such as physics, biology, medicine, or economics, and they constitute an important part of mathematics, naturally putting forth the question of how these functions can be approximated numerically.

Particularly interesting instances are optimization problems where the feasible solutions are a family of convex functions.
For example, let $H^k(\Omega)$ denote the usual Sobolev space of $L^2(\Omega)$ functions having all weak derivatives of order up to $k$ in $L^2(\Omega)$, and suppose $\Omega\subset\RR^d$ is a convex domain.
We may be interested in finding the projection of a given $f\in H^1(\Omega)$ onto the set $\Conv$ of convex functions in $H^1(\Omega)$,
\begin{equation}\label{equ:proj:H1}
   \min_{u\in\Conv} \norm{u - f}_{H^1}.
\end{equation}

Or, given $f \in H^{-1}(\Omega)$, we may be interested in minimizing the Dirichlet functional
\begin{equation}\label{equ:Caffa}
   J_f(u) = \frac{1}{2} \int_\Omega |\grad u(x)|^2\, \dd{x} + \intern{f,u},
\end{equation}
over the set of convex functions $u$ defined in $\Omega$ with $\int_\Omega u\,\dd{x} = 0$.

Often the convexity requirement in applications comes from a reasonable shape assumption on the model, which could be replaced by or added to other shape constraints such as radial symmetry, harmonicity or upper and lower bounds.
This is the case, for instance, of Newton's problem of minimal resistance~\cite{BG97,BFK95,LP99,LP01}.

More surprisingly perhaps, the convexity may be a consequence of the model, as in some mechanism design problems in economics.
For example, Rochet and Chon\'e~\cite{R-C} and Manelli and Vincent~\cite{manelli} (among others) study what we will call the \emph{monopolist problem}, in which the functional to be maximized is the seller's expected revenue,
\begin{equation}\label{equ:monopolist}
   \max_{u\in\Conv}
      \int_Q \bigl(
         \grad u(x) \cdot x - u(x) - c\, \abs{\grad u}^2
         \bigr) \,f(x)\, \dd{x},
\end{equation}
where
$Q = [0,1]^d$ is the $d$-dimensional unit cube,
$c$ is a non-negative constant,
$f$ is a probability density function on $Q$,
and $\Conv$ is the set of convex functions $u$ satisfying
$u(0) = 0$ and $0 \le \grad u \le 1$ (component-wise).
In this problem, the convexity restriction comes from the requirement of incentive compatibility.
We refer the interested reader to~\cite{manelli} and the references therein for further details on the model.

From a theoretical point of view, Carlier and Lachand-Robert~\cite{carlier-cpam} obtained the $C^1$ regularity of a variant of the monopolist problem~\eqref{equ:monopolist}, under some restrictions on the domain $\Omega$ and the density $f$.
They obtained also $C^1$ regularity for convex minimizers of functionals similar to that in~\eqref{equ:Caffa}, with the condition $\int_\Omega u \,\dd{x} = 0$ substituted for $u = u_0$ in $\boundary{\Omega}$.

From a numerical point of view, Carlier, Lachand-Robert and Maury~\cite{carlier} proposed a finite element scheme for minimizers in $H^1(\Omega)$ or $L^2(\Omega)$  of functionals encompassing~\eqref{equ:proj:H1}, \eqref{equ:Caffa} and~\eqref{equ:monopolist}.
In a two dimensional setting, they consider finite element functions which are the interpolants of convex functions, and show that this definition is equivalent to an intrinsic one, stated only in terms of the value of the function at the grid points.
In a way, they consider (weak) second order pure derivatives in every possible direction allowed by the underlying grid.
The problem with this description is that it is non-local, and the number of constraints needed in two dimensions (after pruning) reportedly grows approximately as $N^{1.8}$, where $N$ is the number of nodes in the grid.
Moreover, their approach is very difficult to extend in practice to higher dimensions.

The work of Carlier et al.~\cite{carlier} includes the problem of finding
\[
   \min\int_{\Omega} \abs{u - f}^2\,\dd{x}
   \qquad\text{subject to}\quad
   u\in L^2(\Omega),\ \text{$u$ convex},\ u\le f,
\]
for given $f\in L^2(\Omega)$, i.e., a $L^2$-norm projection, and---as they point out---this problem is equivalent to that of finding the convex envelope $f^{**}$ of $f$.
Thus, minimizing over convex functions and finding the convex hull of the epigraph of a function are two quite related tasks.

Being a central problems in computational geometry, there are a number of well established codes for finding the convex hull of a set of points in $\RR^d$, which are very efficient in low dimensions.
Hence, it is natural to try to use these codes to approximate convex functions, an approach which Lachand-Robert and Oudet~\cite{LO05} applied to several problems.

There is a large literature on convex functions in a continuous setting, well represented by the book by Rockafellar~\cite{rockafellar}.
Also the discrete mathematics community has produced quite a few definitions for convexity of functions defined on lattices (see, e.g., the article by Murota and Shioura~\cite{murota} and the references therein).
But in either case, the definitions are usually of a non-local nature.

One of the main difficulties in obtaining discrete approximations to convex functions in dimensions higher than one, lies in giving a local and finite description of them.
Though this could be done for smooth functions of continuous variables by asking the Hessian matrix to be positive semidefinite at all points, we know of no similar characterizations for finite element functions on meshes.

This article builds on our previous work~\cite{A-M}, where we gave a theoretical framework for approximating convex functions using a finite difference discretization of the Hessian matrix and semidefinite programs.

Let us recall that a semidefinite program is an optimization problem of the form
\begin{equation}\label{equ:psd}
\begin{gathered}
   \min\, c\cdot x\\
   \begin{array}{rc}
   \text{subject to}&\\
   & x_1 A_1 + \dots + x_n A_n - A_0 \succeq 0,\\
   &x\in\RR^n,
   \end{array}
\end{gathered}
\end{equation}
where $c\in\RR^n$, $A_0,A_1,\dots,A_n$ are symmetric $m\times m$ matrices, and $A\succeq 0$ indicates that the symmetric matrix $A$ is positive semidefinite.
By letting the matrices $A_i$ be diagonal, we see that the program~\eqref{equ:psd} is a generalization of linear programming (and includes it strictly).
Thus, in a semidefinite program the constraints can be a mixture of linear inequalities and positive semidefinite requirements.
We refer the reader to the article by Vandenberghe and Boyd~\cite{vandenberghe} for further properties of semidefinite programs.

In this article we carry over the framework in~\cite{A-M} into the approximation with finite elements.
We do so through a weak definition of a finite element (\FE) Hessian and corresponding definition of \FE-convex functions.%
   \footnote{We denote them with the prefix ``\FE'' to distinguish them from other definitions, such as those in~\cite{A-M,carlier,murota}.}
Being the definition very natural and straightforward, the main goal of this article is to provide a solid theoretical foundation of this approach and to illustrate its applicability to a broad range of models.

In contrast to finite differences, it is now very easy to adaptively refine the meshes and reduce drastically the computational effort, especially taking into account the fact that the time needed by the semidefinite programs is more than quadratic on the number of nodes.

Although not linear, our approach seems very natural and has many advantages.
Being of a local nature, the number of constraints grows only linearly with the number of nodes, and it works for any dimension of the underlying space.

The rest of this paper is organized as follows.

In section~\ref{sec:hessian} we introduce the \FE-Hessian and \FE-convex functions, discussing several related issues. We give examples and counter-examples showing how \FE-convex functions relate to usual convex functions and the finite element version given by Carlier et al.~\cite{carlier}.

In section~\ref{sec:limits} we prove the main results of the paper.
We show that, under appropriate assumptions and norm, every convex function, even if not smooth, can be approximated by a sequence of \FE-convex functions, and that the limit of every convergent sequence of \FE-convex functions (with space discretization parameter going to zero) is a convex function.
We also show some compactness results, such as that a (norm) bounded sequence of \FE-convex functions has a convergent subsequence (to a convex function).

In section~\ref{sec:approximating} we show how the previous results may be used to approximate many optimization problems, providing a general framework for the numerical treatment of optimization problems over convex functions, and prove some theoretical results supporting the potential applicability to a broad range of concrete problems.
We do not focus here on a specific problem, and thus our convergence results will not provide convergence rates, since these depend on the regularity properties of the exact solutions, and other features of the particular problem at hand.

In section~\ref{sec:numerical} we discuss the actual numerical implementation, and give concrete examples of the monopolist problem~\eqref{equ:monopolist}, and the Dirichlet integral~\eqref{equ:Caffa}.

We conclude by summing up and commenting on the results we found.

\section{Discrete Hessians and discrete convexity}
\label{sec:hessian}

There are two main issues when defining the set of discrete approximants to be used:

{\enumerate
\item it must be rich enough to approximate every convex function, and

\item it must be not too large, to avoid convergence to non-convex functions.
\endenumerate}

The first point is very natural, and necessary to be able to approximate the solution of the problem.
The second point looks artificial at first sight, but if it did not hold, a sequence of functions in this kind of sets could converge to a non-convex function.

\begin{figure}
\centering
\begin{minipage}[b]{.49\textwidth}\centering
\includegraphics[width=3.5cm]{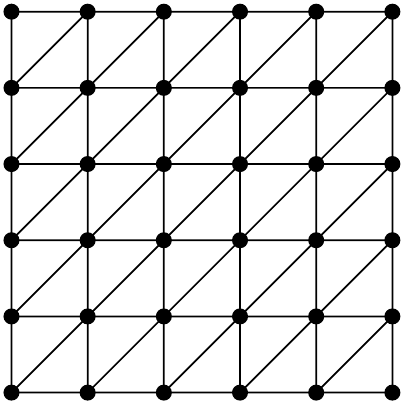}\\
(a) Regular diagonal mesh.
\end{minipage}
\hfill
\begin{minipage}[b]{.49\textwidth}\centering
\includegraphics[width=4cm]{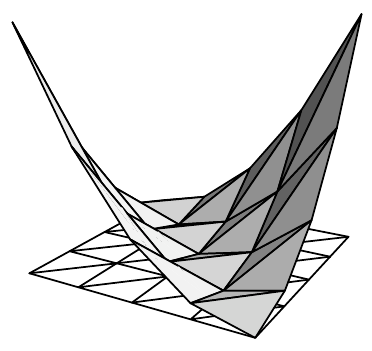}\\
(b) Interpolant of $(x_1 + x_2)^2$.
\end{minipage}
\caption{Interpolant of a convex quadratic function on a regular mesh.}
\label{fig:diagonal}
\end{figure}

On one hand, as noted by P. Chon\'e in his Ph.D. Thesis~\cite{chone}, the affine Lagrange interpolant of a convex function need not be convex.
Consider for instance a regular diagonal mesh as shown to the left in figure~\ref{fig:diagonal}, and suppose $\Omega = (-1,1)\times(-1,1)$.
The interpolant on this mesh of the quadratic convex function $(x_1 + x_2)^2$, shown to the right of that figure, is clearly not convex.

On the other hand, if we consider a convergent sequence of convex piecewise linear functions on a sequence of meshes like those of figure~\ref{fig:diagonal}, with mesh size tending to zero, then the limit will satisfy
\[ \frac{\partial^2 u}{\partial x_1 \partial x_2} \le 0. \]
This is a consequence of proposition~1 in~\cite{carlier}, which was first proved in~\cite{chone}.
It clearly indicates that not all convex functions can be approximated by discrete functions that are convex in the usual sense, i.e., the definition of the discrete approximants to convex functions needs to be more involved.

In order to proceed, we briefly review some concepts and set some notation. 
If $\Omega$ is a bounded open convex set in $\RR^d$ ($d\ge 2$), $u\in C^2(\Omega)$ and $x = (x_1,\dots,x_d)\in\Omega$, the Hessian $Hu(x)\in\RR^{d\times d}$ is defined as the matrix whose $ij$ entry is the second order partial derivative of $u$ in the directions $x_i$ and $x_j$,
\[ {(Hu(x))}_{ij} = \partial_{ij} u (x). \]
As is well known, $u\in C^2(\Omega)$ is convex if and only if $Hu$ is positive semidefinite everywhere in $\Omega$, in symbols
\[ Hu(x)\succeq 0 \quad\text{for all $x\in \Omega$.} \]

When $u$ is not smooth enough, we may nevertheless consider the Hessian in the distribution sense.
In other words, we may define $Hu$ as a matrix of distributions such that for every $\varphi\in C_0^\infty(\Omega)$,
$\intern{Hu,\varphi}$ is the matrix of numbers
\begin{equation}\label{equ:Hess:distr}
   \intern{Hu, \varphi} = \intern{u, H\varphi}.
\end{equation}

In 1971, Dudley~\cite{dudley} proved that if $u$ is a distribution on $\Omega$ such that
\begin{equation}\label{equ:dudley:1}
   \intern{Hu,\varphi} \succeq 0
   \quad\text{for all $\varphi\in C_0^\infty(\Omega)$, $\varphi \ge 0$,}
\end{equation}
then $u$ belongs to the (Lebesgue) equivalence class of a continuous convex function in $\Omega$.
Conversely, if $u$ is a convex continuous function then~\eqref{equ:dudley:1} holds (in the distribution sense).

By allowing some smoothness on $u$, say $u\in H^1(\Omega)$, we may rewrite~\eqref{equ:Hess:distr} and interpret $Hu$ as a matrix of distributions $\begin{bmatrix}(Hu)_{ij}\end{bmatrix}$ satisfying
\begin{equation}\label{equ:Hess:H1}
   \intern{(Hu)_{ij},\varphi}
   = \intern{Hu,\varphi}_{ij}
   = - \int_\Omega \partial_i u(x)\,\partial_j \varphi(x)\,\dd{x},
   \quad\text{for all $\varphi\in C_0^\infty(\Omega)$.}
\end{equation}
In this case, the equality in~\eqref{equ:Hess:H1} also holds for all $\varphi\in H_0^1(\Omega)$, and Dudley's results imply that given $u\in H^1(\Omega)$, $u$ is a continuous convex function in $\Omega$ if and only if
\begin{equation}\label{equ:dudley:2}
   H_v u \succeq 0
   \quad\text{for every $v\in H_0^1(\Omega)$, $v\ge 0$,}
\end{equation}
where for convenience, for $u\in H^1(\Omega)$ and $v\in H_0^1(\Omega)$, we have denoted by $H_v u$ the matrix whose $ij$ entry is
\[
   (H_v u)_{ij} = - \intern{\partial_i u, \partial_j v}.
\]

It is then natural to define a discrete Hessian in the finite element setting along these lines.
To do so, it will be convenient to use two different families of finite element basis functions.
The first one, $\{\phi^h_r\}$ indexed by $r\in\IndV{h}$, will be used for approximations, and the second one, $\{\varphi_s\}$ indexed by $s\in\IndW{h}$, will be used as test functions, and we will assume that
\begin{equation}\label{equ:test:positive}
   \varphi^h_s (x) \ge 0 \quad
   \text{for all $x\in\Omega$ and all $s\in\IndW{h}$.}
\end{equation}

$V_h$ and $W_h$ will denote the (real) linear spaces spanned by $\{\phi^h_r\}$ and $\{\varphi^h_s\}$, respectively, and again for simplicity we will assume $V_h\subset H^1(\Omega)$ and $W_h\subset H^1_0(\Omega)$.
$h$ will denote, as usual, a discretization parameter, equivalent to the maximum diameter of the elements of the underlying grid.

For $u\in V_h$ and each $s\in\IndW{h}$, we define the \emph{\FE-Hessian (of $u$ with respect to $\varphi_s$)}, $H^h_s u$, by
\[ H^h_{s} u = H_{\varphi^h_s} u, \]
and in particular, if $u = \phi^h_r$, we define
\[
   H^h_{rs}
   = H^h_{s}\phi^h_r
   = \begin{bmatrix}
      - \intern{\partial_i \phi^h_r,\partial_j \varphi^h_s}
      \end{bmatrix},
\]
so that
\[ H^h_{s}\Bigl(\sum_r u_r\,\phi^h_r\Bigr) = \sum_r u_r\,H^h_{rs}. \]

We are now in a position to state the following

\begin{definition}\label{def:discrete}
We will say that $u\in V_h$ is \emph{\FE-convex (with respect to $\{\phi^h_r\}$ and $\{\varphi^h_s\}$)} if
\[
   H^h_{s} u
   \succeq 0
   \quad \text{for all $s \in\IndW{h}$.}
\]
\end{definition}

If $u\in V_h$ is convex in the usual sense, and the conditions~\eqref{equ:test:positive} hold, by Dudley's results in the form~\eqref{equ:dudley:2}, $H^h_{s} u\succeq 0$ for all $s\in\IndW{h}$.
Therefore, convex functions in $V_h$ are \FE-convex.

As was shown in the example of figure~\ref{fig:diagonal}, the interpolant of a continuous convex function need not be convex.
However, as we will see next, in that example the interpolant is \FE-convex, and therefore in general \FE-convexity does not imply convexity.

{\example\label{example:diagonal}\upshape
Let us consider a regular diagonal mesh in $\Omega = (0,1)\times (0,1)$, as shown in figure~\ref{fig:diagonal}, let $h$ be the length of the shorter sides of the triangles, and let $V_h$ and $W_h$ consist of piecewise linear functions.

A simple calculation shows that, if $u_h \in V_h$ and $\varphi^h_s\in W_h$ is the function which equals $1$ on the interior node with coordinates $(a,b)$ and vanishes in the other mesh nodes, then
\[
   H^h_{s} u_h =
      \begin{bmatrix} \alpha & \beta \\ \beta & \gamma \end{bmatrix},
\]
where
\begin{align*}
   \alpha &= u_h(a-h,b) + u_h(a+h,b) - 2 u_h(a,b),\\
   \beta &= \frac{1}{2}\,\Bigl(2\,u_h(a,b) + u_h(a-h,b-h) + u_h(a+h,b+h) \\
      &\quad
      -\bigl( u_h(a,b-h) + u_h(a,b+h) + u_h(a-h,b) + u_h(a+h,b)\bigr)
      \Bigr), \\
   \gamma &= u_h(a,b-h) + u_h(a,b+h) - 2 u_h(a,b).
\end{align*}

If $u_h$ is the Lagrange interpolant of the quadratic function $u(x_1,x_2) = (x_1 + x_2 - 1)^2$, another simple calculation shows that
\[ H^h_{s} u_h = \intern{Hu,\phi^h_s} \succeq 0, \]
so that $u_h$ is \FE-convex but it is not convex as we know from figure~\ref{fig:diagonal}.\END
\endexample}

\begin{figure}
\centering
\begin{minipage}[b]{.35\textwidth}\centering
\includegraphics[width=3.5cm]{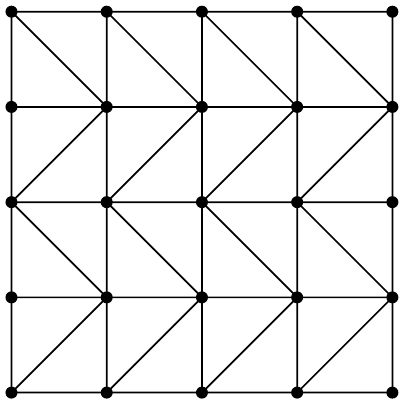}\\
(a) Chevron or alternating.
\end{minipage}\hfill
\begin{minipage}[b]{.3\textwidth}\centering
\includegraphics[width=3.5cm]{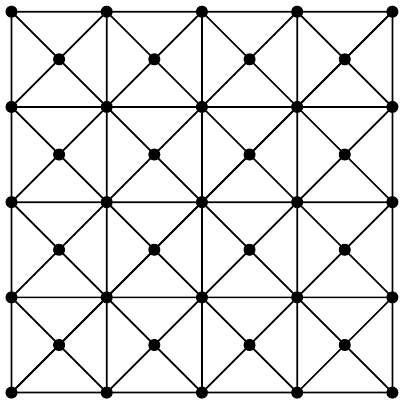}\\
(b) Crisscross.
\end{minipage}\hfill
\begin{minipage}[b]{.3\textwidth}\centering
\includegraphics[width=3.5cm]{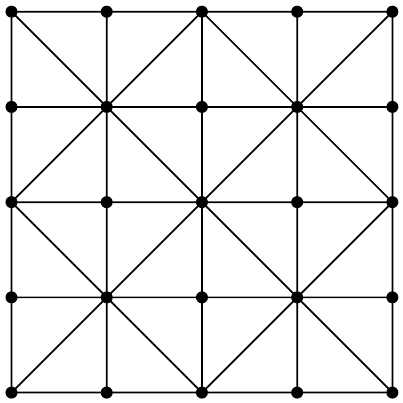}\\
(c) Union Jack
\end{minipage}
\caption{Other common patterns for regular meshes.}
\label{fig:mesh:other}
\end{figure}

It is worth noticing that, in general, it is not true that the interpolant of a convex function is \FE-convex, even for some highly regular meshes.
In order to illustrate this, we have sketched some common patterns of regular meshes in figure~\ref{fig:mesh:other}.
It can be readily seen that the example~\ref{example:diagonal} for the diagonal mesh in figure~\ref{fig:diagonal}(a) can be carried over to the ``chevron'' or ``alternating'' mesh.
However, the behavior is quite different for the ``crisscross'' and ``Union Jack'' patterns.

{\example\label{example:unionjack:1}\upshape
Consider as in the previous example $\Omega = (0,1)\times (0,1)$, but now a regular ``Union Jack'' mesh in $\Omega$ as shown in figure~\ref{fig:mesh:other}(c).
As before, let $h$ be the length of the shorter sides of the triangles, and let $V_h$ and $W_h$ consist of piecewise linear functions.

In this mesh, the nodes inside $\Omega$ can have either $8$ or $4$ neighbors.
If $(a,b)$ is a mesh node having $8$ neighbors and $\varphi_s$ is the corresponding nodal basis function, for $u_h \in V_h$ we now have
\[
   H^h_{s} u_h =
      \begin{bmatrix} \alpha & \beta \\ \beta & \gamma \end{bmatrix},
\]
where
\begin{align*}
   \alpha &= u_h(a-h,b) + u_h(a+h,b) - 2 u_h(a,b),\\
   \beta &= \frac{1}{2}\,\bigl(
      u_h(a+h,b+h) + u_h(a-h,b-h)
      \\&\quad
      - u_h(a+h,b-h) - u_h(a-h,b+h) \bigr), \\
   \gamma &= u_h(a,b-h) + u_h(a,b+h) - 2 u_h(a,b).
\end{align*}

If $u_h$ is the interpolant of the quadratic function $u(x_1,x_2) = (x_1 + x_2)^2$, we see that
\[
   H^h_{s} u_h
   = h^2\,\begin{bmatrix} 2 & 4 \\ 4 & 2 \end{bmatrix} \not\succeq 0,
\]
so that $u_h$ is \emph{not} \FE-convex.

Similar examples can be constructed for the ``crisscross'' meshes since these are essentially $45\degree$ rotations of ``Union Jack'' meshes.\END
\endexample}

We conclude from these examples that our concept of \FE-convex functions neither contains nor is contained in that of~\cite{carlier}.
However, as we show in the next section, these concepts have many common features.

\section{Limits of \FE-convex functions}
\label{sec:limits}

Having defined \FE-convexity through a discrete Hessian, we would like to see how these concepts may be used to approximate convex functions.

The first decision we have to make is in what sense the approximation will be done.
By the very definition of the \FE-Hessian in~\eqref{equ:Hess:H1}, it is natural to consider the approximation in the $H^1$ sense, and this is what we will do, but of course other spaces could be used.
In particular, even when working with approximations in $H^1(\Omega)$, we will make use of the spaces $W^{k,p}(\Omega)$ consisting of the functions having at least $k$ weak derivatives in $L^p(\Omega)$.

On the other hand, it will be convenient to work with a sequence of finer and finer meshes, $\Mesh{h_n}$, with $h_n\downto 0$ as $n\to\infty$.
However, so as not to clutter the notation, we will drop the index $n$, and we will also assume $0 < h < 1$.

We are confronted now with several tasks:

{\enumerate
\item Suppose $(u_h)_h$ is a sequence, $h\downto 0$, of \FE-convex functions $u_h\in V_h$, and suppose that, as $h\downto 0$, $u_h$ has a $H^1(\Omega)$ weak limit $u\in H^1(\Omega)$.
Is $u$ convex?

\item Given a bounded sequence ${(u_h)}_{h}$ with $u_h$ \FE-convex in $V_h$.
Is there a convergent subsequence?

\item Recalling Chon\'e's observations and the example in figure~\ref{fig:diagonal}, can \emph{any} convex function in $H^1(\Omega)$ be approximated as much as desired (in that space) by \FE-convex functions (for appropriate $V_h$ and $W_h$)?
\endenumerate}

The first two issues will be covered by theorem~\ref{theorem:fem:1} and corollary~\ref{coro:fem:1}.
The last issue will be covered by theorem~\ref{theorem:fem:2}, but it is somewhat different in flavor from the previous results, since we will need some properties of the approximating spaces $V_h$.
These will depend on the choice of the finite element families, which can vary widely.

Since it is not our purpose in this paper to present the results under the most general conditions, for simplicity we will restrict our attention to $C^0$ Lagrange elements, hoping that the interested reader will be able to adapt the proof of theorem~\ref{theorem:fem:2} to other families.

We start by observing that if $\Omega$ is any convex domain, then the convex functions in $H^1(\Omega)$ may be approximated by smooth convex functions:

{\lemma\label{lemma:Cinfty}
Suppose $\Omega$ is a bounded convex open subset of $\RR^d$. If $u$ is a convex function in $H^1(\Omega)$, for any $\varepsilon > 0$ there exists a convex $v\in C^\infty(\clausure{\Omega})$ such that
\[ \norm{u - v}_{H^1(\Omega)} < \varepsilon. \]
\endlemma}

{\proof
Let $\varphi\in C_0^\infty$, $\varphi \ge 0$, with support inside $\{x\in\RR^d : \abs{x} < 1\}$ and $\int \varphi\,\dd{x} = 1$, and for $\delta > 0$ consider the mollifier
\[
   u_\delta (x) = u \ast \varphi_\delta (x)
   = \frac{1}{\delta^d}
      \int u(y)\, \varphi(\delta^{-1}\,(x - y))\,\dd{x},
\]
which is well defined in
\[ \Omega_\delta = \{x \in \Omega : \dist(x,\boundary\Omega) > \delta \}, \]
where $\dist(x,A)$ is the distance from $x$ to the set $A$.

As is well known, $u_\delta\in C^\infty(\Omega_\delta)$, and since $u$ is convex and $u_\delta$ is an average with non-negative weights, $u_\delta$ is convex in $\Omega_\delta$.
For $\delta' > 0$ fixed, $u_\delta$ converges to $u$ in $H^1(\Omega_{\delta'})$ as $\delta\downto 0$, and moreover $u_\delta\in C^\infty\bigl(\clausure{\Omega}_{\delta'}\bigr)$ if $\delta < \delta'$.
Thus the result is true for $\Omega_{\delta'}$ for every $\delta' > 0$.

Given its geometry, it is easy to obtain the result for $\Omega$ by using suitable dilations.
For instance, pick $x_0\in\Omega$ and consider for $1 \le \lambda \le 2$, the set $\Omega^\lambda = \{x_0 + \lambda\,(x - x_0): x\in\Omega\}$, and for $u$ defined on $\Omega$ consider $u^\lambda$ defined on $\Omega^\lambda$ by $u^\lambda(x_0 + \lambda\,(x - x_0)) = u(x)$.
Then
$\Omega^\lambda$ is convex,
$\norm{u^\lambda - u}_{H^1(\Omega)} \to 0$ as $\lambda\downto 1$,
and if $u$ is convex in $\Omega$ then $u^\lambda$ is convex in $\Omega^\lambda$.

Thus, fixing first $\lambda > 1$ so that
\(
   \norm{u^\lambda - u}_{H^1(\Omega)} < \varepsilon/2,
\)
then $\delta' > 0$ so that
\(
   {(\Omega^\lambda)}_{\delta'} \supset \Omega,
\)
and finally $\delta > 0$ so that
\(
   \norm{u^\lambda -
      {(u^\lambda)}_\delta}_{H^1\left({(\Omega^\lambda)}_{\delta'}\right)}
      < \varepsilon/2,
\)
we have
\begin{align*}
   \norm{u - {(u^\lambda)}_\delta}_{H^1(\Omega)}
   &\le
      \norm{u - u^\lambda}_{H^1(\Omega)}
      +
      \norm{u^\lambda - {(u^\lambda)}_\delta}_{H^1(\Omega)}
   \\
   &\le \varepsilon/2 + \norm{u^\lambda -
      {(u^\lambda)}_\delta}_{H^1\left({(\Omega^\lambda)}_{\delta'}\right)}
   < \varepsilon.
\end{align*}
\endproof}

From now on we will assume that we have a sequence of meshes $\Mesh{h}$, with $h\downto 0$, each consisting of a family $\Triang{h}$ of non-overlapping closed $d$-dimensional simplices such that for each $h$, $\clausure{\Omega} = \cup_{T\in\Triang{h}} T$.
This implies that $\Omega$ is the interior of a polyhedron (intersection of finitely many half-spaces).
Recalling that we are thinking of a hierarchical sequence of meshes and we require the non-negativity condition~\eqref{equ:test:positive}, we now make some additional assumptions on the meshes and discrete spaces we will work with, following chapter~4 of the book by Brenner and Scott~\cite{B-S}.

{\assumptions\label{assumptions}
We will denote by $\abs{A}$ the Lebesgue measure of $A$, and indicate by $a\simeq b$ that for some positive constants $C_1$ and $C_2$ (independent of $h$ and $u$) we have $C_1 a \le b \le C_2 a$,

\begin{enumerate}
\item
For all $0 < h < h'$, $V_{h'} \subset V_h \subset H^1(\Omega)$ and $W_{h'}\subset W_{h} \subset H^1_0(\Omega)$.

\item\label{assumption:Vh:interp}
There exists a linear operator $\interp{h}$ with values in $V_h$ (the \emph{interpolant}), an integer $m \ge 2$, and a constant $C$ independent of $u$ and $h$, such that
\begin{gather*}
   \norm{u - \interp{h}u}_{H^1(\Omega)}
      \le C h^{m-1}\, \norm{u}_{H^m(\Omega)},
   \\
   \norm{u - \interp{h} u}_{W^{1,\infty}(\Omega)}
      \le C h^{m-1}\, \norm{u}_{W^{m,\infty}(\Omega)}.
\end{gather*}

In particular, $\cup_h V_h$ is dense in $H^1(\Omega)$.

\item\label{assumption:Wh:0} Condition~\eqref{equ:test:positive} holds, i.e.
\[
   \varphi^h_s (x) \ge 0 \quad
   \text{for all $x\in\Omega$, all $h>0$, and all $s\in\IndW{h}$.}
\]

\item\label{assumption:Wh:1}
Given $h_0$ and $s_0 \in \IndW{h_0}$, for every $h < h_0$ there exist coefficients $a_s \ge 0$ such that
\[ \varphi_{s_0} = \sum_{s\in\IndW{h}} a_s\,\varphi_s. \]

\item\label{assumption:Wh:2}
Given $\varphi\in C_0^\infty(\Omega)$, $\varphi\ge 0$, there exists a sequence $(w_h)_h$ converging to $\varphi$ in $H_0^1(\Omega)$ such that
\[ w_h = \sum_{s\in\IndW{h}} a^h_s\,\varphi_s \in W_h, \]
with $a^h_s \ge 0$ for all $h$ in the sequence and $s\in\IndW{h}$.

In particular, $\cup_h W_h$ is dense in $H^1_0(\Omega)$.

\item\label{assumption:Wh:3}
For all $h$,
\begin{align*}
   \abs{T} &\simeq h^d
      &&\text{for all $T\in\Triang{h}$,}
   \\
   \int_\Omega \varphi^h_s\,\dd{x} \simeq h^d
   \quad &\text{and} \quad
   \abs{\grad \varphi^h_s} \le C/h
      &&\text{for all $s\in\IndW{h}$.}
\end{align*}
\end{enumerate}
\endassumptions}

These conditions will be satisfied for quasi-uniform families, taking $C^0$ Lagrange elements with polynomials of degree less than $m$ for the trial space $V_h$, and piecewise linear elements for the test space $W_h$.
The assumptions also hold choosing $W_h$ as the finite element space of continuous piecewise polynomial functions of any fixed degree.
Assumption~\ref{assumptions}.\ref{assumption:Wh:2} is guaranteed because $W_h$ will always contain the piecewise linear finite element functions.
The only detail to take into account is the choice of the basis functions in order to fulfill assumption~\ref{assumptions}.\ref{assumption:Wh:0}.
If the degree is bigger than one, we cannot use as $\varphi^h_s $ the canonical nodal basis functions of $W_h$, because some of them change sign.
The construction is still possible though (see section~\ref{sec:implementation} for details).

The following is one of the main results of the paper.

{\theorem\label{theorem:fem:1}
Let ${(u_h)}_h$ be a sequence converging weakly in $H^1(\Omega)$ to $u$ as $h\downto 0$, such that for each $h$, $u_h\in V_h$ and
\[ H^h_{s} u_h \succeq 0 \quad\text{for all $s\in\IndW{h}$}. \]

Then $u$ is convex.
\endtheorem}

{\proof
By assumption~\ref{assumptions}.\ref{assumption:Wh:1}, for arbitrary but fixed $h_0$ and $s_0\in\IndW{h_0}$, and any $h \le h_0$ we may find coefficients $a^h_s \ge 0$ so that $\varphi_{s_0} = \sum_{s} a^h_s\,\varphi_s$.

Therefore, since $H^h_s u_h\succeq 0$ for every $s\in\IndW{h}$, setting $w = \varphi_{s_0}$ we obtain
\[
   H^h_w u_h = \sum_{s\in\IndW{h}} a^h_s\, H^h_s u_h \succeq 0
   \quad\text{for $h \le h_0$.}
\]

Given that $u_h$ converges weakly in $H^1(\Omega)$ to $u$, we now have
\[ H_w u = \lim_{h\to 0} H^h_w u_h \succeq 0, \]
and finally, by assumption~\ref{assumptions}.\ref{assumption:Wh:2}, $H_\varphi u \succeq 0$ for all non-negative $\varphi\in C_0^\infty(\Omega)$, and the theorem follows by Dudley's results~\eqref{equ:dudley:1}.
\endproof}

Since the unit ball in $H^1(\Omega)$ is weakly compact, we have:

{\corollary\label{coro:fem:1}
Let ${(u_h)}_{h}$ be a bounded sequence in $H^1(\Omega)$  such that for every $h$ the function $u_h \in V_h$ is \FE-convex.
Then there exists a subsequence that converges weakly in $H^1(\Omega)$ to some function $u\in H^1(\Omega)$, and this function is necessarily convex.
\endcorollary}

This theorem and its corollary answer the first two issued raised at the beginning of this section.
In order to proceed further, and answer the last one, we will need the following result by Hoffman and Wielandt~\cite{H-W}:

{\theorem\label{theorem:H-W}
There exists a positive constant $c_d$, depending only on the dimension $d$, such that if $A = [a_{ij}]$ and $B = [b_{ij}]$ are symmetric $d\times d$ matrices, and $\lambda$ and $\mu$ are their minimum eigenvalues, then
\[ \abs{\lambda - \mu} \le c_d\,\max_{ij}\,\abs{a_{ij} - b_{ij}}. \]
\endtheorem}

The following is the second main result and responds the third issue raised at the beginning of this section.

{\theorem\label{theorem:fem:2}
If $m > 2$, given $u\in H^1(\Omega)$, $u$ convex and $\varepsilon > 0$, there exist $h > 0$ and $u_h\in V_{h}$ such that
\[ \norm{u - u_h}_{H^1(\Omega)} < \varepsilon \]
and
\[ H^h_{s}(u_h) \succeq 0 \quad\text{for all $s\in\IndW{h}$}. \]
\endtheorem}

{\proof
By lemma~\ref{lemma:Cinfty}, it will be enough to assume that $u$ is a $C^\infty (\clausure{\Omega})$ convex function.

In the sequel we will denote by $C$, $C'$ or $C''$, positive constants which may vary from one occurrence to another, even in the same line, which may depend on $u$ (which we consider fixed from now on), $\Omega$, the dimension $d$ and the regularity degree $m$, but are independent of $h$.
For instance, we write
\[ \norm{u}_{H^m(\Omega)} = C. \]

Let us consider the auxiliary function
\[ g(x) = \frac{1}{2}\,\abs{x}^2, \]
which is a convex $C^\infty(\RR^d)$ function.
The regular Hessian of $g$, $Hg$, is
\[ Hg = I_d = \text{identity matrix in $\RR^{d\times d}$}, \]
and therefore, for any $w\in H^1_0(\Omega)$,
\[
   H_w g = \intern{Hg, w}
   = \int_{\Omega} I_d\, w\,\dd{x}
   = \bigl(\int_\Omega w\,\dd{x}\bigr)\,I_d,
\]
and in particular,
\begin{equation}\label{equ:res:1}
   H^h_{s} g =  \Bigl(\int_\Omega \varphi^h_s\,\dd{x}\Bigr)\,I_d,
   \quad\text{for all $h$ and all $s\in\IndW{h}$.}
\end{equation}

For $\delta$ and $h$ positive and small, let
\[ v = u + \delta g \qquad\text{and}\qquad u_h = \interp{h} v, \]
where $\interp{h}$ denotes the interpolant considered in assumption~\ref{assumptions}.\ref{assumption:Vh:interp}.
We notice that since the third derivatives of $g$ vanish, $h$ is bounded and $m\ge 2$,
\begin{gather*}
   \norm{g}_{H^m(\Omega)} = \norm{g}_{H^2(\Omega)} = C,
   \\
   \norm{v}_{H^m(\Omega)}
      \le \norm{u}_{H^m(\Omega)} + \delta\,\norm{g}_{H^m(\Omega)}
      \le C,
   \\
   \norm{\interp{h}g}_{H^1(\Omega)}
      \le C\, \norm{g}_{H^m(\Omega)} = C,
\end{gather*}
and
\begin{align*}
   \norm{u - u_h}_{H^1(\Omega)}
   &
   = \norm{u - \interp{h}u - \delta\,\interp{h}g}_{H^1(\Omega)}
   \le \norm{u - \interp{h}u}_{H^1(\Omega)}
      + \norm{\delta\,\interp{h}g}_{H^1(\Omega)}
   \\
   &
   \le C h^{m-1}\, \norm{u}_{H^m(\Omega)} + C\,\delta
   \le C\, (h + \delta).
\end{align*}

Thus, if $h + \delta$ is small enough (depending on $\varepsilon$), we obtain the first inequality of the theorem.

In order to see that $H^h_{s} u_h \succeq 0$, we first look at the eigenvalues of $H^h_{s} v$ for $s\in\IndW{h}$.
If $\zeta\in\RR^d$, using that $u$ is convex and smooth, equation~\eqref{equ:res:1}, and the bounds in assumption~\ref{assumptions}.\ref{assumption:Wh:3}, we obtain
\begin{multline*}
   \bigl(H^h_{s} v\, \zeta\bigr)\cdot \zeta
   = \bigl(\intern{Hu, \varphi^h_s}\, \zeta\bigr)\cdot \zeta
      + \delta\, \bigl(\intern{Hg, \varphi^h_s}\, \zeta\bigr)\cdot \zeta
   \\
   \ge 0 + \delta\, \abs{\zeta}^2 \int_{\Omega} \varphi^h_s\,\dd{x}
   \ge C \delta\,\abs{\zeta}^2\,h^d,
\end{multline*}
and therefore the eigenvalues of $H^h_{s} v$ are bounded below by
\begin{equation}\label{equ:res:3}
   C\delta h^d.
\end{equation}

In order to compare the entries of $H^h_s v$ and $H^h_s u_h$, we use assumptions~\ref{assumptions}.\ref{assumption:Wh:3} and~\ref{assumptions}.\ref{assumption:Vh:interp} to obtain
\begin{equation}\label{equ:fem:2}
   \begin{aligned}
   \bigl\lvert{(H^h_{s} u_h - H^h_{s} v)}_{ij}\bigr\rvert
   &
   = \bigl\lvert{(H^h_{s} \interp{h} v - H^h_{s} v)}_{ij}\bigr\rvert
   = \Bigl\lvert -\int_\Omega
      \partial_{i} (\interp{h} v - v)\, \partial_j \varphi^h_s \,\dd{x}
      \Bigr\rvert
   \\
   &
   \le \norm{\interp{h} v - v}_{W^{1,\infty}(\Omega)}\,
      \norm{\varphi^h_s}_{W^{1,1}(\Omega)}
   \\
   &
   \le C' h^{m-1} \norm{v}_{W^{m,\infty}(\Omega)}\, h^{d-1}
   = C' h^{m + d - 2}.
   \end{aligned}
\end{equation}
Thus we may use theorem~\ref{theorem:H-W} and the bounds~\eqref{equ:res:3} and~\eqref{equ:fem:2} to obtain that the eigenvalues of $H^h_s u_h$ are bounded below by
\[ C \delta h^d - C' h^{m + d - 2} = C'' h^d\,(\delta - h^{m-2}). \]

The theorem follows now by taking $\delta \ge C h^{m-2}$ if $m > 2$, and at the same time $\delta + h$ small.
\endproof}

The condition $m > 2$ in theorem~\ref{theorem:fem:2} implies that the functions in $V_h$ will have to be piecewise polynomials of degree at least $2$, meaning that the result may not hold for linear finite elements.
This does not seem quite satisfactory, and we need to elaborate on the necessity of this condition.

As we have seen in example~\ref{example:diagonal}, for meshes such as those of figure~\ref{fig:diagonal}(a) or figure~\ref{fig:mesh:other}(a), if we use piecewise linear functions for the space $V_h$ ($m = 2$), the discrete Hessian becomes a finite difference scheme (except for a factor of $h^2$) which is exact for quadratic functions.
The results presented in~\cite{A-M} can be adapted to show that in this case theorem~\ref{theorem:fem:2} also holds for $m=2$.

On the other hand, as we have seen in example~\ref{example:unionjack:1}, the discrete Hessian with piecewise linear functions for very regular meshes such as those in figure~\ref{fig:mesh:other}(b) or (c) is not exact for quadratic functions, which means that we may not be able to get good approximations.
In the following example we report numerical evidence supporting the necessity of assuming $m>2$ in theorem~\ref{theorem:fem:2}.

\begin{figure}
\centering
\includegraphics[width=.3\textwidth]{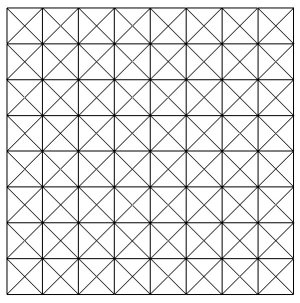}\hfill
\includegraphics[width=.3\textwidth]{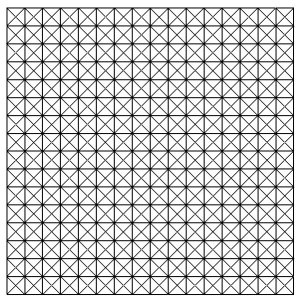}\hfill
\includegraphics[width=.3\textwidth]{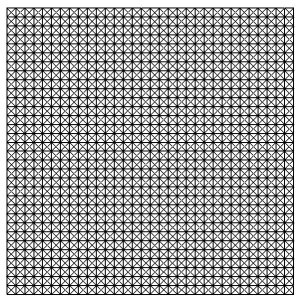}
\\[4pt]
\includegraphics[width=.3\textwidth]{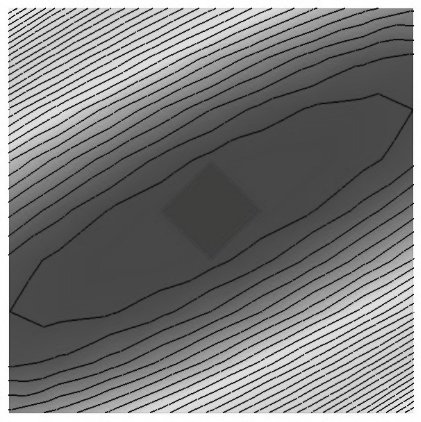}\hfill
\includegraphics[width=.3\textwidth]{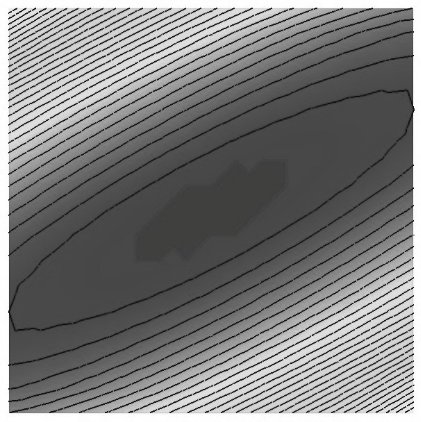}\hfill
\includegraphics[width=.3\textwidth]{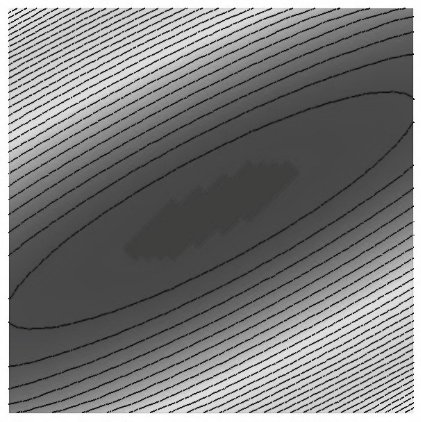}
\caption{%
Example showing that convergence cannot be guaranteed when using piecewise linear functions.
Meshes (top), and $L^2$-projection into \FE-convex piecewise linear functions (bottom).
The level curves of the projected function---which is already convex and should be the limit of the discrete ones---are straight lines, which are not reproduced by the approximants.
The meshes have 256, 1024, and 4096 elements.}
\label{fig:noconv}
\end{figure}

{\example\label{example:unionjack:2}\upshape
We compute the $L^2( (0,1)\times(0,1) )$-projection of the smooth convex function
\[ u(x_1,x_2) = (x_2 - 0.5 x_1 -0.25 )^2, \]
onto the set of continuous piecewise linear functions that are \FE-convex over criss-cross meshes, as those obtained using longest-edge or newest-vertex bisection.
Since $u$ is convex, the projections $u_h$ should converge to $u$ as $h\to 0$, but this is not the case in this example.

In figure~\ref{fig:noconv} we can observe a sequence of meshes and the level curves of the $L^2$-projection of $u$ into the set of \FE-convex piecewise linear functions on each mesh.
The level curves of the exact function $u$, which is convex, are straight lines parallel to the line $x_2 = 0.5 x_1$.
Nevertheless, the level curves of the approximants converge to ellipses which do not straighten up by refinement.
This is a clear indication that $(u_h)$ does not converge to $u$ as $h \to 0$. We remark that the same behavior is observed when projecting in $H^1$, $L^\infty$ and $L^1$, and even when imposing boundary values.\END
\endexample}

Summing up: although there is some sort of super-convergence for some meshes, for general meshes---even highly regular---\FE-convex piecewise \emph{linear} functions may not suffice to approximate convex functions.

\section{Approximating functionals}
\label{sec:approximating}

We are now in position to apply finite element approximations to a wide class of optimization problems on convex functions.

Let us describe this technique by assuming, for instance, that the functional
\[ J(v) = \int_\Omega F(x,v(x),\grad v(x)) \,\dd{x} \]
is defined and continuous on $H^1(\Omega)$, and we are interested in the optimization problem
\begin{equation}\label{prob:PJ}
   \inf\, \{ J(v): v\in \Conv\},
\end{equation}
where $\Conv$ is a family of convex functions, $\Conv\subset H^1(\Omega)$.

Using theorem~\ref{theorem:fem:2}, it may be not too difficult to define for each $h > 0$
a family $\Conv_h\subset V_h$, and a functional $J_h$ defined on $\Conv_h$, such that:
{\enumerate
\item $H_s^h v_h(x) \succeq 0$ for all $v_h\in\Conv_h$ and $s \in \IndW{h}$,

\item for any $v\in \Conv$ and any $\varepsilon > 0$, there exists $h > 0$ and $v_h\in \Conv_h$ such that $\abs{J_h(v_h) - J(v)} < \varepsilon$.
\endenumerate}

Under the previous conditions, it is easy to prove that (cf.~\cite{A-M})
\begin{equation}\label{equ:fund}
   \inf\, \{ J(v): v\in \Conv\} = \lim_{h\to 0}\,\inf\, \{ J_h(v_h) : v_h\in \Conv_h \}.
\end{equation}

As a concrete example, suppose that $\Conv$ is the set of all convex functions in $H^1(\Omega)$ with a given mean value, or some prescribed boundary values, $f\in H^1(\Omega)$, and the continuous problem consists in finding $u\in\Conv$ such that
\[ \norm{u - f}_{H^1} = \min_{v\in\Conv} \norm{v - f}_{H^1}, \]
i.e., minimizing $ J(v) := {\norm{v - f}}_{H^1}^2$ over $\Conv$.

In order to compute an approximation of $u$ we may thus consider $\Conv_h$ as the set of discrete functions $v_h\in V_h$ with $H_s^h v_h\succeq 0$ for all $s\in\IndW{h}$, which also satisfy the integral or boundary constraints.

Assuming exact integration we can set $J_h(v_h) = J(v_h)$, or otherwise, $J_h(v_h)$ may result from some fixed quadrature rule on the elements of the mesh.
In both cases it is easy to see that the previous assumptions hold, and thus the discrete minimizers $u_h \in \Conv_h$ provide a convergent (sub)sequence to the exact solution $u$.

\section{Numerical Experiments}
\label{sec:numerical}
\subsection{Implementation issues}\label{sec:implementation}

The numerical examples were implemented using the finite element toolbox \textsf{ALBERTA}~\cite{alberta} for assembling the optimization problem, and CSDP~\cite{csdp1} for solving the corresponding semidefinite programs.
The experiments were run on a desktop PC with a 2.8~GHz Intel Pentium IV processor and 2GB of RAM.

In our experiments we used Lagrange finite elements of polynomial degree $2$ for both $V_h$ and $W_h$, over simplicial meshes, but the right choice might depend on the precise problem at hand.

Regarding the implementation in \textsf{ALBERTA}, we had to introduce some modified basis functions when using quadratic finite elements as test functions in $W_h$.
The canonical nodal basis functions associated with the vertices of the elements change sign, and the theory requires that the test functions be non-negative.
To do this, we considered the usual piecewise \emph{linear} nodal basis functions for the vertices, whereas for the nodes that correspond to the midpoints of the edges we used the usual quadratic bubbles, which are obtained as the product of the two linear basis functions that correspond to the vertices of the edge.

In the initial experiments we observed some oscillations at the boundary, but they ceased to appear when we incorporated into the test space the basis functions corresponding to the boundary nodes, enlarging $W_h$ so that it is no longer a subset of $H_0^1(\Omega)$.
In this case, formula~\eqref{equ:Hess:H1} was transformed into
\begin{equation}\label{equ:Hess:H1:bdry}
   \intern{(Hu)_{ij},\varphi}
   = \intern{Hu,\varphi}_{ij}
   = - \int_\Omega \partial_i u(x)\,\partial_j \varphi(x)\,\dd{x}
      + \int_{\partial\Omega} \partial_i u(x) \, \varphi(x) \nu_j\,
      \dd{S},
\end{equation}
where $\nu_j$ denotes the $j$-th component of the outward unit normal to $\partial\Omega$.
This slight modification still leads to the same theoretical results of the previous sections.
We decided to present them assuming zero boundary values for the test functions, in order to keep the presentation clearer.

\subsection{Statement of the discrete problems}

In the examples that follow, we always considered the minimization of functionals of the form
\[
   J(u)
   = \int_\Omega \bigl(
      \alpha\, \abs{\grad (u - v_1)}^2
      + \beta\, \abs{u - v_2}^2
      + \gamma \cdot \grad u + f u
      \bigr) \,\dd{x},
\]
where $\alpha$, $\beta$, $\gamma$, $v_1$, $v_2$, and $f$ are given functions on $\Omega$.
Appropriate choices of these functions lead to functionals whose minima are the $L^2(\Omega)$- or the $H^1(\Omega)$-projection of a function, or the solutions to problems~\eqref{equ:Caffa}, or~\eqref{equ:monopolist}.
The approximate functional $J_h(v_h)$ results from applying some fixed quadrature rule on the elements of the mesh.

In order to model this problem as a semidefinite program on each given mesh, we used a fixed quadrature rule (exact for polynomials of degree $\le 4$) over the elements of the mesh and approximated the functional by
\begin{multline*}
   J_h(u) = \sum_{i=1}^N w_i \bigg[
      \alpha(x_i) \sum_{j=1}^d
      \big(\frac{\partial u}{\partial x_j}(x_i)
      - \frac{\partial v_1}{\partial x_j}(x_i)\big)^2 \\
      + \beta(x_i) ( u(x_i) - v_2(x_i) )^2
      + \gamma(x_i) \cdot \grad u(x_i) + f(x_i) u(x_i)
      \bigg],
\end{multline*}
where $x_i$, $w_i$, $i=1,2,\dots,N$ are the quadrature points and weights, respectively.
The minimization of $J_h$ was then modeled by adding $(d+1)N$ auxiliary variables $t_{ij}$, $s_i$, $i=1,2,\dots,N$, $j=1,2,\dots,d$, as
\[
\text{minimize } \sum_{i=1}^N w_i \bigg[ \alpha(x_i) \sum_{j=1}^d t_{i,j}
                            + \beta(x_i) s_i
                            + \gamma(x_i) \cdot  \grad u(x_i) + f(x_i) u(x_i)
\bigg]
\]
subject to
\begin{align*}
\big(\frac{\partial u}{\partial x_j}(x_i) - \frac{\partial v_1}{\partial x_j}(x_i)\big)^2 &\le t_{i,j}, \\
( u(x_i) - v_2(x_i) )^2 &\le s_i
\end{align*}
for $i=1,2,\dots,N$, and $j=1,2,\dots,d$, plus the convexity constraints.
In turn, the constraints involving squares are modeled, respectively, by
\[
\left[\begin{matrix} 1 & \frac{\partial u}{\partial x_j}(x_i) - \frac{\partial v_1}{\partial x_j}(x_i) \\
      \frac{\partial u}{\partial x_j}(x_i) - \frac{\partial v_1}{\partial x_j}(x_i) & t_{i,j}
    \end{matrix}
    \right]  \succeq 0
\]
and
\[
\left[\begin{matrix} 1 & u(x_i) - v_2(x_i) \\
      u(x_i) - v_2(x_i) & s_i
    \end{matrix}
    \right]  \succeq 0.
\]

\subsection{Adaptivity}

In order to take full advantage of the flexibility of finite elements, we included some adaptivity into our algorithms, which was implemented as a loop of the form~\cite{mns}
\[
\textsc{solve } \to\textsc{ estimate }\to\textsc{ mark }
\to\textsc{ refine}.
\]
The step \textsc{solve} consisted in solving the resulting semi-definite programs using CSDP.
Having computed the discrete solution $u_h$, the step \textsc{estimate} consisted in estimating the error distribution over the triangulation $\Triang{h}$ in the following way: we defined for each $T\in\Triang{h}$ the quantity $\eta_T = h_T^{1/2} \| [ \grad u_h ] \|_{L^2(\partial T)}$, where $[ \grad u_h ]_{|S}$ denotes the jump of $\grad u_h$ over the inter-element sides $S$ and is defined as zero at boundary sides.
This quantity $\eta_T$ is the dominating part of the residual-type a posteriori error estimator for Poisson's problem, and we used it as a heuristic indicator of the error.
Further studies are necessary in order to develop rigorous upper and lower bounds for the error in this type of problems, an open question which is out of the scope of this article; we introduced a heuristic error estimator here just to show the power of finite elements and the great improvement in performance that adaptivity can provide.
The step \textsc{mark}, consisted in marking (selecting for refinement) all the elements whose indicators satisfied $\eta_T \ge 0.7 \eta_{\text{max}}$, where $\eta_{\text{max}} := \max_{T \in \Triang{h}} \eta_T$.
The step \textsc{refine} was implemented using the standard routines of \textsf{ALBERTA}, which perform newest-vertex bisection, guaranteeing a uniform shape-regularity constant.

\subsection{Examples}

\begin{figure}
\centering
\includegraphics[width=.3\textwidth]{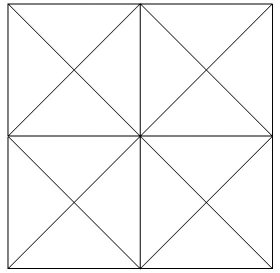}\hfill
\includegraphics[width=.3\textwidth]{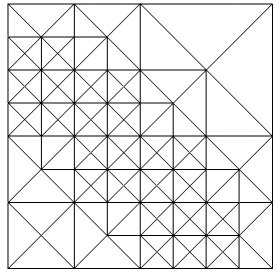}\hfill
\includegraphics[width=.3\textwidth]{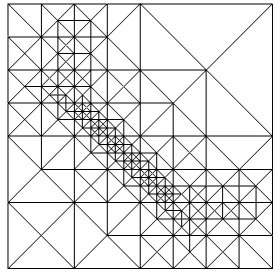}\\[4pt]
\includegraphics[width=.3\textwidth]{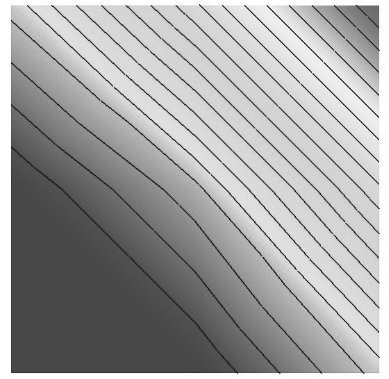}\hfill
\includegraphics[width=.3\textwidth]{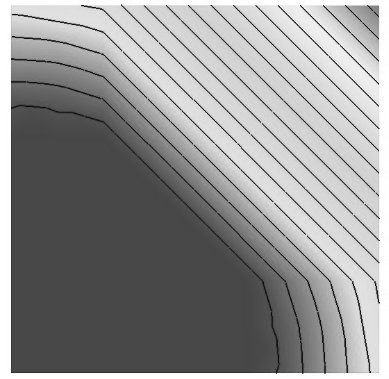}\hfill
\includegraphics[width=.3\textwidth]{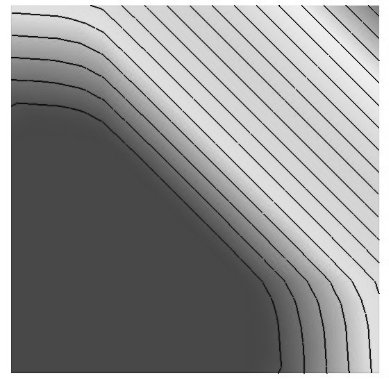}
\caption{%
Contour lines and meshes obtained with the adaptive version of the algorithm.
We show the meshes of iterations 0, 2, and 4, with 16 (41), 140 (299), and 388 (799), elements (DOFs), respectively.}
\label{fig:monopolist}
\end{figure}

\begin{figure}
\centering
\includegraphics[width=.3\textwidth]{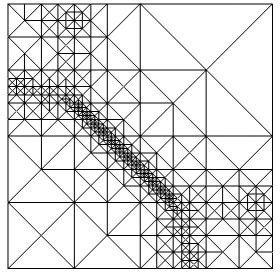}\hfill
\includegraphics[width=.3\textwidth]{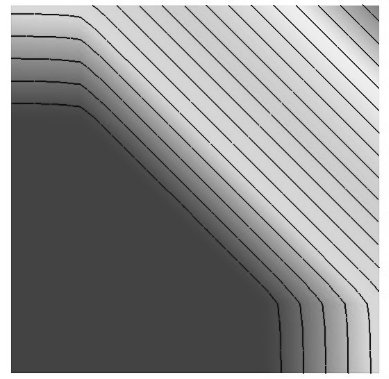}\hfill
\includegraphics[width=.3\textwidth]{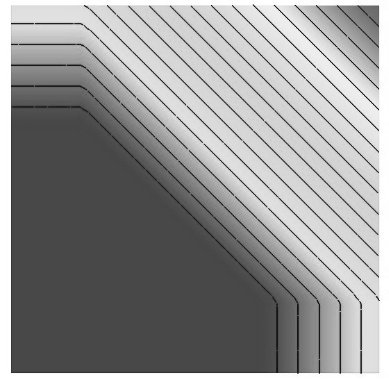}
\caption{%
Final mesh, approximate and exact solution, after 6 iterations.
The mesh has 802 elements and 1641 degrees of freedom were used with quadratic elements.
The code CSDP took 211 seconds to minimize the functional, and the $L^2$-error between the approximate and exact solution is $0.00192$.
To obtain a similar error with uniform meshes, more than 4096 elements and 8321 DOFs are needed, which forces CSDP to work more than 7000 seconds (see table~\ref{tab:monopolist-unif}).}
\label{fig:monopolist-final}
\end{figure}

{\example\label{example:monopolist}\upshape
In this example we apply our algorithm to solve the monopolist problem~\eqref{equ:monopolist}, for $d=2$, $f\equiv 1$ and $c=0$.
In this case the exact solution is known to be
\[ u(x_1,x_2) = \max\,\{ 0, x_1-a, x_2-a, x_1+x_2-b \}, \]
where $a = 2/3$ and $b = (4-\sqrt2)/3$, and allows us to compute the true error.
The method is applied using quadratic elements both in the trial and test spaces.
In figure~\ref{fig:monopolist} we show a sequence of solutions using adaptive meshes, and in table~\ref{tab:monopolist} we can observe the error. In Figure~\ref{fig:monopolist-final} we show the final mesh, approximate and exact solution, after 6 iterations.
In order to illustrate the performance of the adaptive method, we also include in table~\ref{tab:monopolist-unif} the errors and CPU-times obtained with uniform meshes. The reported CPU-times correspond to the time taken by CSDP to find the minimum of the functional on the given mesh. To be fair in the comparison, we should look at the cumulative sum of the CPU-times, since the whole adaptive process is necessary in order to arrive at the graded meshes.\END
\endexample}

\begin{table}
\centering
\begin{tabular}{|r|r|r|r|r|}
\hline
\strutmio[-4pt]{16pt}%
Elements & DOFs & CPU-time & $L^2$-error & $L^\infty$-error \\
\hline
\strutmio{12pt}%
    16  &    41  &    0.00 & 0.06371 & 0.11168 \\
    54  &   123  &    1.00 & 0.01597 & 0.06214 \\
   140  &   299  &   11.00 & 0.00935 & 0.02589 \\
   248  &   519  &   25.00 & 0.01075 & 0.01553 \\
   388  &   799  &   59.00 & 0.01023 & 0.01633 \\
   542  &  1117  &   94.00 & 0.01085 & 0.01529 \\
   802  &  1641  &  211.00 & 0.00192 & 0.00846 \\
\hline
\end{tabular}\\[10pt]
(a)\\[12pt]
\begin{tabular}{|r|r|r|r|r|}
\hline
\strutmio[-4pt]{16pt}%
Elements & DOFs & CPU-time & $L^2$-error & $L^\infty$-error \\
\hline
\strutmio{12pt}%
    16  &    41 &     1 & 0.06371 & 0.11168 \\
    64  &   145 &     2 & 0.01644 & 0.06250 \\
   256  &   545 &    22 & 0.01202 & 0.02347 \\
  1024  &  2113 &   343 & 0.01091 & 0.01780 \\
  4096  &  8321 &  7211 & 0.01270 & 0.03508 \\
\hline
\end{tabular}\\[10pt]
(b)
\caption{Errors and CPU-time (in seconds) for the adaptive (a)
and uniform-refinement (b) solutions of the monopolist problem.}
\label{tab:monopolist}
\label{tab:monopolist-unif}
\end{table}

{\example\label{example:Caffa}\upshape
In this example we apply our algorithm to minimize the functional defined in~\eqref{equ:Caffa} over the set of convex functions in $\left\{ v\in H^1(\Omega):\int_\Omega v = 0\right\}$ (cf.~\cite{carlier-cpam}).
We consider $\Omega = \left\{ (x_1,x_2) \in \RR^2 : x_1^2 + x_2^2 < 1\right\}$ and

\medskip

\noindent
\begin{minipage}[c]{.75\textwidth}
\[
   f(x_1,x_2) = \begin{cases}
      1 & \text{if } x_1^2 + (x_2+1)^2 \le 1/4, \\
      -1 & \text{if } x_1^2 + (x_2-1)^2 \le 1/4, \\
      0 & \text{otherwise.}
      \end{cases}
\]
\end{minipage}
\hfill
\begin{minipage}[c]{.2\textwidth}
\includegraphics[width=.7\textwidth]{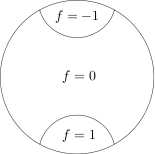}
\end{minipage}

\medskip

In figure~\ref{fig:Caffa} we show the outcome of our method using adaptivity.
We can observe that the solution \emph{tries} to satisfy $\Delta u = f$ in places where $f \ge 0$, and it continues to be convex outside that region, minimizing $|\Delta u - f|$.
In the lower part of the domain (around the point $(0,-1)$), the solution $u$ satisfies $\Delta u = f$ and on the boundary, the natural homogeneous Neumann boundary conditions ${\partial u}/{\partial \nu} = 0$ are satisfied.
This can be seen by the fact that the level curves are perpendicular to the boundary in that part of $\partial \Omega$.
In the upper part of the domain the solution is just linear (ruled surface), which is a consequence of the fact that the Laplacian of $u$ has to be non-negative and as close to $-1$ as possible, keeping $u$ convex over the whole domain.

\begin{figure}
\centering
\includegraphics[width=.3\textwidth]{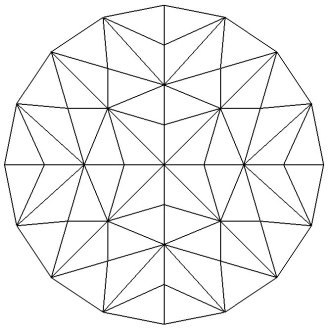}\hfill
\includegraphics[width=.3\textwidth]{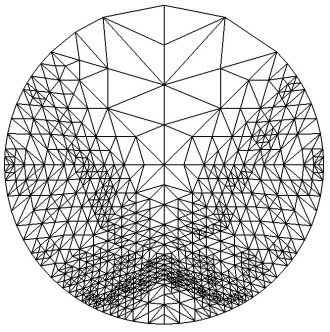}\hfill
\includegraphics[width=.3\textwidth]{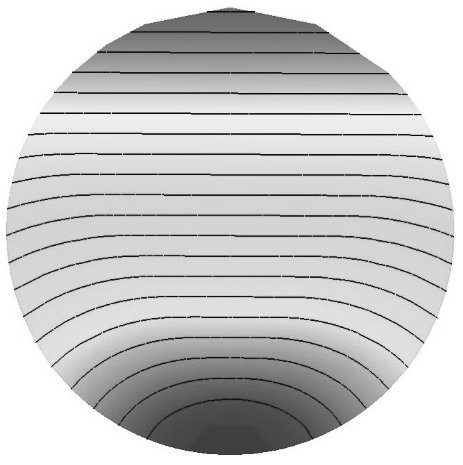}\\[4pt]
\includegraphics[width=.31\textwidth]{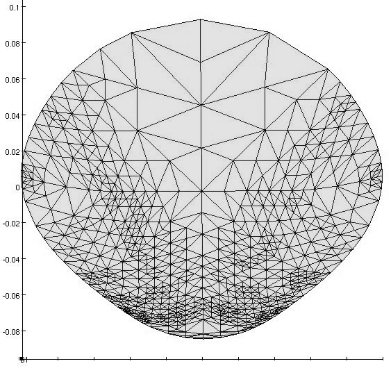}\hfill
\includegraphics[width=.31\textwidth]{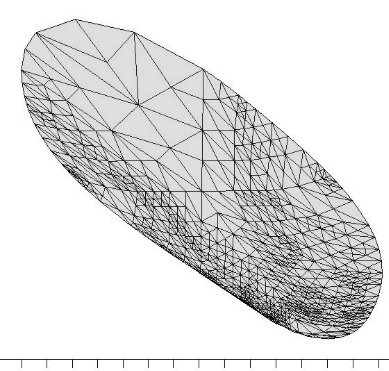}\hfill
\includegraphics[width=.31\textwidth]{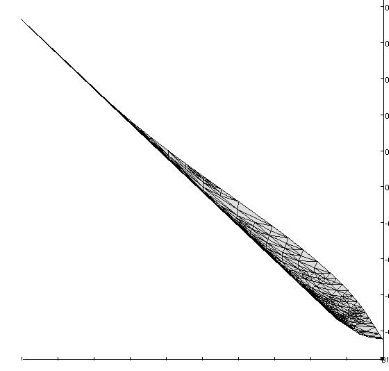}
\caption{%
Minimizer of the functional~\eqref{equ:Caffa}.
Top: initial and final mesh, and contour curves of the solution.
Bottom: surface plot of the solution as viewed from $(0,-10,0)$ (left), $(-10,-10,0)$ (middle), and $(-10, 0, 0)$ (right).}
\label{fig:Caffa}
\end{figure}

The adaptive method correctly captures the region where $u$ is flat, representing the solution with a minimal number of elements.\END
\endexample}

\section{Conclusions}
\label{sec:conclusions}

We have proposed a novel way of imposing convexity on finite element functions, and proved that this new definition solves the two issues necessary for the approximation of optimization problems over convex functions:

\begin{itemize}
\item Every convex function can be approximated; and
\item If a sequence of \FE-convex functions is convergent, then the limit is convex.
\end{itemize}

Our definition takes advantage of the existence of efficient codes for solving semidefinite programs, and uses a new definition of discrete Hessians, which is based upon a weak Hessian for $H^1$ functions.
One interesting aspect is that the definition is intrinsic, i.e., it only uses the values of the discrete functions at the nodes, and \emph{local}, leading to a set of constraints with cardinality of order $O(N)$, $N$ being the number of vertices or nodes of the mesh.
Furthermore, it is very simple to program in any space dimension.

Another interesting---and puzzling---issue is the fact that, in general, except for some particular very regular meshes, the discrete functions need to have an approximation order higher than the one provided by linears.
Our proof requires this assumption, and we found numerical evidence that this is necessary, but a better explanation/understanding of this issue is still pending.

Numerical experiments show a competitive performance, specially through the use of adaptivity, which, in turn, is easy to implement for finite elements.
Our preliminary computations using a heuristic error indicator are promising, but a lot needs to be done in this direction, namely, to find a posteriori error indicators which are reliable and efficient, and once this is established, to prove convergence, and optimality.
These are difficult open questions that will be subject of future research.



\end{document}